
\input amstex

\documentstyle{amsppt}

\loadbold

\magnification=\magstep1

\pageheight{9.0truein}
\pagewidth{6.5truein}


\topmatter 

\title LENGTHS OF FINITE DIMENSIONAL REPRESENTATIONS 
OF PBW ALGEBRAS \endtitle

\date September 2003 \enddate 

\author D. Constantine and M. Darnall \endauthor

\address Eastern Nazarene College,
23 E. Elm Ave.,
Quincy, MA 02170 \endaddress 

\address Humbolt State University,
1 Harpst St.,
Arcata, CA 95521 \endaddress

\abstract
Let $\Sigma$ be a set of $n\times n$ matrices with entries from a
field, for $n > 1$, and let $c(\Sigma )$ be the maximum length of
products in $\Sigma$ necessary to linearly span the algebra it
generates.  Bounds for $c(\Sigma )$ have been given by Paz and
Pappacena, and Paz conjectures a bound of $2n-2$ for any set of
matrices.  In this paper we present a proof of Paz's conjecture for
sets of matrices obeying a modified Poincar\'e-Birkhoff-Witt (PBW)
property, applicable to finite dimensional representations of Lie
algebras and quantum groups.  A representation of the quantum plane
establishes the sharpness of this bound, and we prove a bound of
$2n-3$ for sets of matrices which do not generate the full algebra of
all $n\times n$ matrices.  This bound of $2n-3$ also holds for
representations of Lie algebras, although we do not know whether it is
sharp in this case.\endabstract

\thanks This research was begun during the 2003 Temple Mathematics Research
Experience for Undergraduates, supported by Temple University and NSF REU Site
Grant DMS-0138991. The authors were participants in this program. \endthanks

\endtopmatter

\document

\head 1. INTRODUCTION \endhead 

For a fixed integer $n > 1,$ let $\Sigma=\{X_1, \ldots , X_t\}$ be a set of $n\times n$ matrices over an
arbitrary field {\bf k}, and let $\Sigma^m$ be the set of products of length $m$ in
the $X_i,$ where $\Sigma^0$ is defined as the identity. Let $L_i$ be the
linear space spanned by $\Sigma^0\cup\Sigma^1\cup \ldots \cup\Sigma^i,$ and
denote the dimension of this space by $r_i.$
Next, let $L_*$ be the linear space spanned by products of any length,
and let $r_*$ denote its dimension.  Finally, let $c(\Sigma)=min\{i:r_i=r_*\}.$

In [3], Paz proved that $c(\Sigma)\leq \lceil (n^2+2)/3\rceil$, and
Pappacena gave lower bounds [2]. Paz conjectured a bound of $2n-2$ and
suggested a lemma which, if proved, would prove the conjecture.  We
prove this lemma (listed as Proposition 2.5 below) for matrices
satisfying the following property: every product
$u=X_{i_1}X_{i_2}\cdots X_{i_l}$ in the matrices $X_1, \ldots ,X_t$
can be written, modulo $L_{l-1},$ in the form
$$\sum_{j_1+j_2\cdots +j_t=l} c_{(j_1, \ldots ,
j_t)}X_t^{j_t}X_{t-1}^{j_{t-1}}\cdots X_1^{j_1},$$ 
with $c_{(j_1,
\ldots , j_t)}=0$ whenever $X_t^{j_t}X_{t-1}^{j_{t-1}}\cdots
X_1^{j_1}<u$ in the lexicographical ordering.

This modified PBW property allows any $l$-length matrix product $u$ to
be written as a linear combination of ordered products of length $l$,
modulo products (not necessarily ordered) of lesser length; here by
{\it ordered product} we mean a product in which $X_iX_j$ never
appears for $i<j.$ Our condition is, in fact, looser than what is
generally found in homomorphic images of algebras satisfying the PBW
property.  Sets of matrices obeying this property include finite
dimensional representations of Lie algebras and quantum groups (see,
e.g., [1] for further details).

The $2n-2$ bound is in fact sharp for the class of matrices satisfying the
above property, as an example
using the quantum plane illustrates, but lower bounds for certain other cases
can be obtained.  In particular, knowledge about $r_*$ allows the constraints
of Paz's suggested lemma to be tightened, resulting in lower bounds on
$c(\Sigma).$ As an example we provide a proof of a $2n-3$ bound when the
{\bf k}-algebra generated by $\Sigma$ is not equal to $M_n(\text{{\bf k}}),$ the full
algebra of all $n\times n$ matrices over the field {\bf k}.

This research was begun during the summer 2003 Research Experience for Undergraduates
program supervised by E. S. Letzter at Temple University.  The authors are
greatly indebted to Dr. Letzter for his guidance during the program and his
help in writing this paper. 

\head  2. PROOF OF THE MAIN THEOREM \endhead 

We begin with some notation and preliminaries necessary to our proof.  We
then prove Lemmas 2.3 and 2.4.  Together, these preliminary lemmas establish
Paz's suggested lemma, listed below as Proposition 2.5.  We then proceed to
prove our main theorem.

\subhead 2.1 Notation\endsubhead

(i) We write the matrix product $X_{i_1} \cdots
X_{i_k}$ as the word $i_1 \ldots i_k.$ An {\it $l$-subword} of a word $u$ is any set
of $l$ consecutive letters in $u$.  

(ii) We say two words are {\it formally equivalent} if their $i^{th}$ letters match for all $i;$
otherwise we say they are {\it formally distinct}.  We call a subword consisting
of one repeated letter, such as $111\ldots 1,$ {\it formally
constant}. 

(iii) If $u_i, u_j, u_k, \ldots , u_p$ are $m$-length products, $u_i \propto (u_j, u_k,
\ldots , u_p )$ means $u_i$ is a linear combination of $u_j, u_k, \ldots , u_p$ modulo
$L_{m-1}.$ 

(iv)  We call a word {\it reducible} if it can be written as a
linear combination of words of lesser length.

\subhead 2.2 Preliminaries\endsubhead

(i) Any word that can be written as a linear combination, modulo words of lesser length, of other words that are
all reducible is itself reducible.  Thus, we will
examine only ordered products, and the modified PBW property satisfied by
the matrices we are considering guarantees that our results will carry over to all
matrix products.

(ii) If any word contains a subword of $n$ or more of the same
letter this word will be reducible by the Cayley-Hamilton Theorem.

(iii) Any word $u$ can be treated as a base-$(t+1)$ number, and
this number will be unique to $u.$  We denote this number by $\bar u.$  The
numerical ordering on $\bar u$ coincides exactly with the lexicographical
ordering on matrix products.

\proclaim{Lemma 2.3} For any positive integers $k, m$ and $N$, if $r_k-r_{k-1}\leq N$ then any word of length $m$ with
more than $N$ formally distinct $k$-subwords can be written as a linear combination,
modulo $L_{m-1}$, of words each having at most $N$ formally distinct
$k$-subwords. \endproclaim

\demo{Proof} Given such a word $u$ of length $m$, let $u_1, u_2, \ldots , u_{s}$ for
$s>N$ be $u$'s formally distinct $k$-subwords, numbered such that $\bar u_1 < \bar u_2 <
\cdots < \bar u_{s}.$ Since there are more subwords than $r_k-r_{k-1},$ these
subwords must be linearly dependent, modulo $L_{k-1}.$ Therefore there exists a
minimum $i$ such that $u_i \propto (u_{k_1}, u_{k_2}, \ldots , u_{k_p} )$ with $k_1,
k_2, \ldots , k_p > i.$ (Note that if $i=s$ then $u_s$ is equal, modulo
$L_{k-1},$ to zero, and $u$ is trivially reducible.)  We form the new
words $u^{(1)}, u^{(2)}, \ldots , u^{(p)}$ from $u$ by replacing
$u_i$ with $u_{k_l}$ to form $u^{(l)}.$ We see that $u \propto (u^{(1)},
u^{(2)}, \ldots , u^{(p)})$ and that $\bar u^{(l)} > \bar u$ for all $l.$

We can apply the above process to each of the words $u^{(l)}$ as long as they
have more than $N$ formally distinct subwords.  Since this process continually
increases the numerical value of these words and there are finitely many
words of length $m,$ we will eventually write $u$ as a linear combination, modulo $L_{m-1},$
of words with at most $N$ formally distinct $k$-subwords.

In addition, because rewriting words in ordered form via the modified PBW property also
continually increases their numeric value, we can, at each step in the above
process, put all our words in ordered form.  Therefore, when working with sets
of matrices obeying the modified PBW property we can write any ordered word
$u$ as a linear combination, modulo $L_{m-1},$ of {\it ordered} words with at
most $N$ formally distinct $k$-subwords.\qed\enddemo

\proclaim{Lemma 2.4} For a positive integer $k\leq 2n-2,$ set 
$$N = \cases  k & \text{ for $1\leq k \leq n-1$ } \\
             2n-k-2 & \text{ for $n\leq k \leq 2n-2.$ } \endcases $$
Any ordered word of length $2n-1$
not reducible by the Cayley-Hamilton Theorem contains at least $N+2$ formally
distinct $k$-subwords. \endproclaim

\demo{Proof}  Let
$u$ be an ordered word of length $2n-1$ that is not reducible by the Cayley-Hamilton Theorem.

{\sl Case I.} Suppose $1 \leq k \leq n-1.$ Recall $N=k.$

{\sl Subcase i.} The longest formally constant subword in $u$ has length greater than or equal to $k$.

Specifically, call the longest formally constant subword $w$ and say it has length $j\geq k$.  Since our word is
not reducible by Cayley-Hamilton, $j< n.$  Since $k\leq j<n$ and $u$
has length $2n-1$ there will be at least $k+1=N+1$ $k$-subwords overlapping but
not contained in $w$. Examining
these
$k$-subwords, we
see that they will be
formally distinct since
each features the transition between $w$ and the
surrounding letters in a different spot.   Thus, including one of the formally
constant $k$-subwords found within $w$, we conclude that $u$
contains at least $N+2$ formally distinct $k$-subwords.

{\sl Subcase ii.} The longest formally constant subword in $u$ has length less than $k.$

In this case no two $k$-subwords will be formally equivalent since none will
be constant.  Since $k< n,$ $u$ has at least $n+1$ $k$-subwords.  Since
$N=k < n,$ $u$ contains at least $N+2$ formally distinct $k$-subwords. 

{\sl Case II.} Now suppose $n\leq k \leq 2n-2.$ Recall $N=2n-2-k.$ Since $k\geq n,$ no
two $k$-subwords can be formally equivalent or else $u$ will be reducible by
Cayley-Hamilton since $u$ will contain a formally constant subword of length greater
than $n$.  There are $2n-k$ $k$-subwords in total, so $u$ contains $N+2$
formally distinct $k$-subwords. \qed\enddemo

We now prove the lemma suggested by Paz.

\proclaim{Proposition 2.5} Let $m=2n-1.$ If for some positive integer $k\leq2n-2$
the corresponding condition from among
$$r_k-r_{k-1}\leq k \hbox{  for  } 1 \leq k \leq n-1 $$
$$r_k-r_{k-1}\leq 2n-k-2 \hbox{  for  } n \leq k \leq 2n-2$$ 
holds then $c(\Sigma)\leq m-1.$ \endproclaim

\demo{Proof}  Lemmas 2.3 and 2.4 establish this proposition in the following
manner.  Suppose one of the above conditions holds; say it is the
condition for $k^*$ and let $N^*$ correspond to $k^*$ as described
in the statement of Lemma 2.4.  Consider a word
$u$ of length $m=2n-1.$ We will show that $u$ is reducible, giving us that
$c(\Sigma)\leq m-1 = 2n-2.$ 

If $u$ is reducible by Cayley-Hamilton, we are done.  
If $u$ is not reducible by
Cayley-Hamilton then Lemma 2.4 implies that it has more than
$N^*$ formally distinct $k^*$-subwords.  Lemma 2.3 then implies that $u$ is linearly
dependent, modulo $L_{m-1},$ on ordered words which do not have more than $N^*$ distinct
$k^*$-subwords.  Finally, the contrapositive of Lemma 2.4 implies that these words are reducible by
Cayley-Hamilton.  Thus $u$ is reducible, and $c(\Sigma)\leq m-1.$ \qed\enddemo

We now prove our main theorem.

\proclaim{Theorem 2.6} Let $\Sigma=\{ X_1,\ldots X_t\} $ be a set of $n\times
n$ matrices satisfying the following property: any product $X_{i_1}\cdots X_{i_l}$ can be
written, modulo $L_{l-1},$ in the form 
$$\sum_{j_1+\cdots +j_t=l} c_{(j_1, \ldots j_t)}X_t^{j_t}\cdots
X_1^{j_1},  \text{ with $c_{(j_1, \ldots j_t)}=0$ whenever $\overline{X_t^{j_t}\cdots
X_1^{j_1}}  <\overline{u}$. }$$
Then $c(\Sigma)\leq 2n-2.$ \endproclaim

\demo{Proof (following Paz {\rm [3])}}  If $c(\Sigma)=m\geq2n-1,$ none of the
conditions of Proposition 2.5 can hold.  Thus, if $c(\Sigma)\geq2n-1$, then $r_0=1, r_1-r_0\geq 2, r_2-r_1\geq 3,
\ldots , r_{n-1}-r_{n-2}\geq n, r_n-r_{n-1}\geq n-1, \ldots , r_{2n-2}-r_{2n-3}\geq 1.$
Then we have $r_{2n-2}\geq 1+2+ \cdots + n-1 + n + n-1 + \cdots +1 = 2(n(n-1))/2 + n
=n^2 \geq r_*.$  This, however, contradicts $c(\Sigma)\geq 2n-1,$ for if
$r_{2n-2}$ is already greater than or equal to the dimension of all of $L_*,
r_{2n-1}$ can be no larger than $r_{2n-2}.$ \qed\enddemo

\head 3. SHARPNESS OF THE BOUND \endhead

The bound of $2n-2$ is sharp for the general set of matrices described above
as the following example from the quantum plane shows.

Consider complex $n\times n$ matrices $X$ and $Y$ satisfying $XY=qYX$,
where $q=e^{2\pi i/n}$, such that the algebra generated by $X$ and $Y$
is all of $M_n(\Bbb{C} ).$ Because $X^n$ and $Y^n$ are reducible by
the Cayley-Hamilton Theorem, the set $P = \{X^iY^j|0\leq i, j
\leq n-1\}$ must span all of $L_*$.  Since
$M_n(\Bbb{C} )$ has dimension $n^2$ and $P$ contains $n^2$ matrices,
$P$ is in fact a basis.  Thus the $(2n-2)$-length product
$X^{n-1}Y^{n-1}$ is linearly independent from any products of lesser
length, giving us that $c(\Sigma)=2n-2$ for such a set of matrices.

It remains only to show that such matrices do
indeed exist.  We leave it to the reader to verify that the following matrices
satisfy the above conditions. 
$$X= \bmatrix 1 & & & \\
               & q & & \\
               & & \ddots \\
               & & & q^{n-1} \endbmatrix, \quad
Y= \bmatrix  & & & 1 \\
              1 \\
              & \ddots \\
              & & 1  \endbmatrix.$$

Lower bounds for certain sets of matrices can be obtained.  Paz's
suggested lemma is set up to deal with sets of matrices for which $r_*$ could
be as great as $n^2.$ With more information about the dimension of $L_*$ for a
given set of matrices, the conditions of the lemma can be tightened, resulting
in lower bounds for $c(\Sigma),$ as the following shows.

We prove a slightly more restrictive form of Proposition 2.5, which we then
use to prove Theorem 3.2.

\proclaim{Proposition 3.1} Let $m=2n-2.$ If for some $k\leq2n-3$ the corresponding
condition from among
$$r_k-r_{k-1}\leq k \hbox{  for  } 1 \leq k \leq n-1 $$
$$r_k-r_{k-1}\leq 2n-k-2 \hbox{  for  } n \leq k \leq 2n-3$$ 
holds then $c(\Sigma)\leq m-1.$ \endproclaim

\demo{Proof} Lemma 2.4 tells us that in an ordered word of length $2n-1$ there are at least
$N+2$ formally distinct
$k$-subwords.  Since decreasing to length $2n-2$
eliminates at most one of these $k$-subwords, there will still be at least
$N+1.$  Since we only need more than $N,$ the proof of Proposition 3.1 then follows directly from the proof of
Proposition 2.5.\qed\enddemo

\proclaim{Theorem 3.2} Let $\Sigma$ be as before, with the added restriction
that it does not generate all of $M_n(\text{{\bf k}}).$  Then $c(\Sigma)\leq 2n-3.$ \endproclaim

\demo{Proof} Now we proceed as before.  If $c(\Sigma)=m\geq 2n-2,$ then none of the above
conditions can hold.  This implies $r_0=1, r_1-r_0\geq 2, r_2-r_1\geq 3,
\ldots , r_{n-1}-r_{n-2}\geq n, r_n-r_{n-1}\geq n-1, \ldots , r_{2n-3}-r_{2n-4}\geq 2.$
Then we have $r_{2n-3}\geq 1+2+ \cdots + n-1 + n + n-1 + \cdots +2 =
2(n(n-1))/2 + n - 1
=n^2-1.$ Because of the restriction placed on the algebra generated by $\Sigma,$ $n^2-1\geq r_*.$  As before,
this contradicts $c(\Sigma)\geq 2n-2.$ \qed\enddemo

For representations of Lie algebras $c(\Sigma)$ is bounded by $2n-3$ as
well.  No Lie algebra consisting of  two matrices generates all of $M_n,$ and for a
Lie algebra of three or more matrices to do so those three
matrices, together with the identity, must be linearly independent, implying that $r_1-r_0\geq 3.$ This allows
a proof similar to that given for Theorem 3.2 since now the sum of the
$r_{i+1}-r_i$ terms will be greater than or equal to $n^2.$
We leave as an open question whether the bound of $2n-3$ is sharp for
representations of Lie algebras.  We have looked for an example achieving this
bound but have been unable to find one.


\bigskip

\Refs

\ref \no 1 \by A. Joseph \book Quantum Groups and Their Primitive Ideals
\bookinfo
Ergebnisse der Mathematik und ihrer Grenzgebiete 3 \vol 29 \publ
Springer-Verlag \publaddr Berlin \yr 1995 \endref

\ref \no 2 \by C. J. Pappacena \paper
An Upper Bound for the Length of a Finite-Dimensional
Algebra \jour J. Algebra \vol 197 \yr 1997 \pages 535-545 \endref

\ref \no 3 \by A. Paz \paper An Application of the Cayley-Hamilton Theorem to
Matrix Polynomials in Several Variables \jour J. Lin$.$ Mult$.$ Algebra \vol
15 \yr 1984 \pages 161-170 \endref

\endRefs
\enddocument